\documentclass[11pt]{article}

\usepackage{amsfonts}
\usepackage{amsmath}
\usepackage{amsmath,amssymb,amsthm}

\usepackage{color}
\usepackage{xcolor}
\usepackage{enumerate}
\usepackage{ifpdf}
\usepackage{algorithmic}
\usepackage{float}

\newtheorem{thm}{Theorem}[section]
\newtheorem{lem}[thm]{Lemma}

\newtheorem{rem}[thm]{Remark}

\newtheorem{defn}[thm]{Definition}
\newtheorem{assum}[thm]{Assumption}

\DeclareSymbolFont{symbolsC}{U}{txsyc}{m}{n}
\DeclareMathSymbol{\coloneqq}{\mathrel}{symbolsC}{66}

\ifpdf
  \usepackage[pdftex,unicode]{hyperref}
  \usepackage[pdftex]{graphicx}
\else

    \usepackage[dvipdfm
         ]{graphicx}
    \usepackage[dvipdfm,unicode]{hyperref}
  \AtBeginDvi{} 
\fi

\hypersetup{backref,
            CJKbookmarks=true,
            bookmarks=true,
            pdfstartview=FitH,
            colorlinks=true,
            bookmarksnumbered=true,
            bookmarksopen=true,
            bookmarksopenlevel=4,
            linkcolor=blue,%
            citecolor=blue,%
            plainpages=false,%
            unicode }

\usepackage{listings}
\newcommand{\ba}{\begin{array}}
\newcommand{\ea}{\end{array}}
\newcommand{\bea}{\begin{eqnarray}}
\newcommand{\eea}{\end{eqnarray}}
\newcommand{\beaa}{\begin{eqnarray*}}
\newcommand{\eeaa}{\end{eqnarray*}}
\def\q{\quad}
\def\qq{\qquad}

\def\si{\sigma}
\def\ms{\medskip}
\def\bs{\bigskip}
\def\reff#1{{\rm(\ref{#1})}}
\def\a{\alpha}
\def\b{\beta}
\def\g{\gamma}
\def\d{\delta}
\def\e{\varepsilon}

\def\si{\sigma}
\def\t{\tau}
\def\f{\varphi}

\def\o{\omega}

\def\cd{\cdot}
\def\cds{\cdots}
\def\pa{\partial}
\def\cd{\cdot}
\def\cds{\cdots}
\def\pa{\partial}

\def\olu{\overline{u}}
\def\ulu{\underline{u}}
\def\tow{(t,\o)}

\def \ol{\overline}
\def \ul{\underline}

\def\L{\Lambda}
\def \zo {\mathbf{0}}
\def \dbf{{\mathbf{d}}}
\def\1{{\bf 1}}
\def\ch{\textsc{h}}
\def \Usup{\overline{\cU}}
\def \Usub{\underline{\cU}}

\def\O{\Omega}
\newcommand{\la}{\langle}
\newcommand{\ra}{\rangle}

\def\dbE{\mathbb{E}}
\def\dbF{\mathbb{F}}

\def\dbL{\mathbb{L}}

\def\dbP{\mathbb{P}}
\def\dbR{\mathbb{R}}
\def\dbS{\mathbb{S}}
\def\dbT{\mathbb{T}}

\def\cA{{\cal A}}

\def\cD{{\cal D}}
\def\cE{{\cal E}}
\def\cF{{\cal F}}

\def\cH{{\cal H}}

\def\cL{{\cal L}}

\def\cP{{\cal P}}

\def\cT{{\cal T}}
\def\cU{{\cal U}}

\def\L{\Lambda}

\begin{document}
\title{\bf Monotone Schemes for Fully Nonlinear  Parabolic Path Dependent PDEs}
\author{  Jianfeng {\sc Zhang}\footnote{University of Southern California, Department of Mathematics, jianfenz@usc.edu.}
       \and Jia {\sc Zhuo}\footnote{University of Southern California, Department of Mathematics, jiazhuo@usc.edu.} 
       }
\date{}
\maketitle
\begin{abstract} In this paper we extend the results of the seminal work Barles and Souganidis \cite{BS} to  path dependent case. Based on the viscosity theory of path dependent PDEs, developed by Ekren, Keller, Touzi and Zhang \cite{EKTZ} and Ekren, Touzi and Zhang \cite{ETZ0, ETZ1, ETZ2}, we show that a monotone scheme converges to the unique viscosity solution of the (fully nonlinear) parabolic path dependent PDE. An example of such monotone scheme is proposed. Moreover, in the case that the solution is smooth enough, we obtain the rate of convergence of our scheme. 
\end{abstract}


\ms



\section{Introduction}
\label{sect-Introduction}
\setcounter{equation}{0}
In this paper we aim to numerically solve the following fully nonlinear Path Dependent PDE (PPDE, for short) with terminal condition $u(T,\o) = g(\o)$:
\bea
\label{PPDE}
 \cL u(t,\o):=-\partial_t u (t,\o) - G(t,\omega, u ,\partial_\o u, \partial_{\o\o}^2 u)=0,~0\leq t<T.
\eea
Here $\o$ is a continuous path on $[0, T]$, and $G$ is increasing in $\partial_{\o\o}^2 u$ and thus the PPDE is parabolic. Such PPDE provides a convenient tool for non-Markovian models, especially in stochastic control/game with diffusion control and financial models with volatility uncertainty. Its typical examples include: martingales as path dependent heat equations, Backward SDEs of Pardoux and Peng \cite{PP} as semilinear PPDEs, and $G$-martingales of Peng \cite{Peng-G} and Second Order Backward SDEs of Soner, Touzi and Zhang \cite{STZ} as path dependent HJB equations. The notion of PPDE was proposed by Peng \cite{Peng-ICM}. Based on the functional It\^{o} calculus, initiated by Dupire \cite{Dupire} and further developed by Cont and Fournie \cite{CF}, Ekren, Keller, Touzi and Zhang \cite{EKTZ} and Ekren, Touzi and Zhang \cite{ETZ0, ETZ1, ETZ2} developed a viscosity theory for  PPDEs.

In the Markovian case, namely $u(t,\o) = v(t, \o_t)$, $g(\o)= f(\o_T)$, and $G(t,\o, y,z,\g) = F(t,\o_t,y,z,\g)$ for some deterministic functions $v, f, F$, the PPDE \reff{PPDE} becomes a standard PDE with terminal condition $v(T,x) = f(x)$:
\bea
\label{PDE}
 \dbL v(t,x):=-\partial_t v (t,x) - F(t, x, v , D v, D^2 v)=0,~0\leq t<T.
\eea
In their seminal work Barles and Souganidis \cite{BS} proposed some time discretization scheme for the above PDE and showed that, under certain conditions, the discretized approximation converges to the unique viscosity solution of the PDE. Their key assumption is the monotonicity of the scheme, see Theorem \ref{thm-PDE} (ii) below, which can roughly be viewed as the comparison principle for the discretized PDE. This work has been extended by many authors, either by improving the error analysis including the rate of convergence, or by proposing specific algorithms which indeed satisfy the required conditions, see e.g. \cite{BJ, BZ, FTW, GZZ, Krylov, Tan, Tan2}, to mention a few. 

Our goal of this paper is to extend the work \cite{BS} to PPDE \reff{PPDE}. Notice that the viscosity solution in \cite{EKTZ, ETZ0, ETZ1, ETZ2} is defined through some optimal stopping problem under nonlinear expectation, which is different from the standard viscosity theory for PDEs. 
Consequently, our notion of monotonicity for the scheme also involves the nonlinear expectation, see \reff{mon} below. This requires some technical estimates for the hitting time involved in the theory. Then, following the arguments in \cite{BS} we show that our monotone scheme converges to the unique viscosity solution of the PPDE. 

We next propose a specific scheme which satisfies all the conditions and thus indeed converges. Moreover,  when the PPDE has smooth enough classical solution, we obtain the rate of convergence of our scheme.

In the semilinear case, there have been many works on numerical methods for the associated backward SDEs, see e.g. \cite{BDM, BL, HTT, HNS, MPST, PX,  Zhang}.  In particular, \cite{HTT} used the arguments for viscosity theory of PPDEs. Moreover, \cite{Tan} studied certain numerical approximation for path dependent HJB equations, in the language of second order BSDEs. However, we should point out that most of these works are mainly theoretical studies and are not feasible, especially in high dimensions. Efficient numerical algorithms for path dependent PDEs, including  the implementation of our discretization scheme in the present paper, remains a challenging problem and we shall explore further in our future research.

The rest of the paper is organized as follows. In \S \ref{sect-Preliminary} we introduce path dependent PDEs and its viscosity solutions, as well as monotone schemes for (standard) PDEs. In \S \ref{sect-mon} we prove the main theorem, namely the convergence of monotone schemes. In \S \ref{sect-scheme} we propose a scheme which satisfies all the desired conditions. Finally in \S \ref{sect-rate} we obtain the rate of convergence of our scheme in the case that the solution is smooth enough. 

\section{Preliminaries}
\label{sect-Preliminary}
\setcounter{equation}{0}

\subsection{Path dependent PDEs and  viscosity solutions}
  In this subsection, we recall the setup and the notations of \cite{ETZ0, ETZ1, ETZ2}. 

\subsubsection{The canonical setting}

Let $\O:= \big\{\o\in C([0,T], \dbR^d): \o_0={\bf 0}\big\}$, the set of continuous paths starting from the origin, $B$ the canonical process, $\dbF$ the natural  filtration generated by $B$, 
$\dbP_0$ the Wiener measure,  and $\L :=[0,T]\times \O$. Here and in the sequel, for notational simplicity,  we use ${\bf 0}$ to denote vectors or matrices with appropriate dimensions whose components are all equal to $0$. Let $\dbS^d$ denote the set of $d\times d$ symmetric matrices, and
 \beaa
 &x \cd x' := \sum_{i=1}^d x_i x'_i 
 ~~\mbox{for any}~~x, x' \in \dbR^d,
 ~~\g : \g' := \mbox{Trace}[\g\g']
 ~~\mbox{for any}~~\g, \g'\in \dbS^d.
 &
 \eeaa
We define a semi-norm on $\O$ and a pseudometric on $\L$ as follows: for any $(t, \o), ( t', \o') \in\L$,
\bea\label{rho}
 \|\o\|_{t} 
 := 
 \sup_{0\le s\le t} |\o_s|,
 \q  
 \dbf\big((t, \o),( t', \o')\big) 
 := 
\sqrt{ |t-t'|} + \big\|\o_{.\wedge t} - \o'_{.\wedge t'}\big\|_T.
 \eea
Then $(\O, \|\cd\|_{T})$ is a Banach space and $(\L, \dbf)$ is a complete pseudometric space.  

\begin{rem}
\label{rem-dinfinity}
{\rm In \cite{ETZ0, ETZ1, ETZ2}, following \cite{Dupire} we used pseudometric:
\beaa
 \dbf_\infty\big((t, \o),( t', \o')\big) 
 := 
 |t-t'| + \big\|\o_{.\wedge t} - \o'_{.\wedge t'}\big\|_T.
\eeaa
Clearly $\dbf$ and $\dbf_\infty$ induce the same topology, and all the results in \cite{ETZ0, ETZ1, ETZ2} still hold true under $\dbf$. However, when we consider the regularity of viscosity solutions, see \reff{reg} below, it is more natural to use $\dbf$. Indeed, since $B$ is typically a semimartingale, for $t<t'$ we see that $\sqrt{t'-t}$ and $\|B^t\|_{t'}$ are roughly in the same order.
\qed}
\end{rem}
We shall denote by $\dbL^0(\cF_T)$ and $\dbL^0(\L)$ the collection of all $\cF_T$-measurable random variables and $\dbF$-progressively measurable processes, respectively.  In particular, for any $u\in \dbL^0(\dbF)$, the progressive measurability implies that $u(t,\o) = u(t, \o_{\cd\wedge t})$. Let $C^0(\L)$ (resp. $UC(\L)$) be the subset of $\dbL^0(\L)$  whose elements are continuous (resp. uniformly continuous) in $(t,\o)$ under $\dbf$. The corresponding subsets of bounded processes are denoted as $C^0_b(\L)$ and  $UC_b(\L)$. Finally, $\dbL^0(\L, \dbR^d)$ denote the space of $\dbR^d$-valued processes with entries in $\dbL^0(\L)$, and we define similar notations for the spaces $C^0$, $C^0_b$, $UC$, and $UC_b$.

We denote by $\cT$ the set of $\dbF$-stopping times, and $\cH\subset\cT$ the subset of those hitting times $\ch$ of the form
\bea
\label{ch}
\ch := \inf\{t: B_t \notin O\} \wedge t_0 = \inf\{ t: d(\o_t, O^c) = 0\} \wedge t_0,
\eea
for some $0< t_0\le T$, and some open and convex set $O  \subset \dbR^d$ containing ${\bf 0}$.

For all $L>0$, let $\cP_L$ denote the set of probability measures $\dbP$ on $\O$ such that there exist 
$\a^\dbP\in \dbL^0(\L, \dbR^d)$, ${\bf 0} \le \b^\dbP \in \dbL^0(\L, \dbS^d)$ satisfying
\bea
\label{cPL}
\left.\ba{c}
|\a^\dbP|\le L, \q  |\b^\dbP|\le \sqrt{2L}, \\
M^\dbP_t:= B_t - \int_0^t \a^\dbP_sds ~\mbox{is a $\dbP$-martingale with}~ d\la M^\dbP\ra_t = {1\over 2}(\b^\dbP_t)^2dt, ~\mbox{$\dbP$-a.s.}
\ea\right.
\eea
and we define $\cP_\infty := \bigcup_{L>0} \cP_L$. We note that, when $\b^\dbP>{\bf 0}$, the second line above is equivalent to the existence of  a $d$-dimensional $\dbP$-Brownian motion $W^\dbP$ satisfying:
\bea
\label{cPL0}
 d B_t =  \a^\dbP_t dt + \b^\dbP_t dW^\dbP_t,\q \dbP\mbox{-a.s.}
\eea

We define the path derivatives via the functional It\^{o} formula.

\begin{defn}
\label{defn-spaceC12}  We say $u\in C^{1,2}(\L)$ if $u\in C^0(\L)$ and there exist $\pa_t u \in C^0(\L)$, $\pa_\o u \in C^0(\L, \dbR^d)$, $\pa^2_{\o\o} u\in C^0(\L, \dbS^d)$ such that, for any $\dbP\in \cP_\infty$, $u$ is a local  $\dbP$-semimartingale and it holds:
\bea
\label{Ito}
d u = \pa_t u dt+ \pa_\o u \cd d B_t + \frac12 \pa^2_{\o\o} u : d \la B\ra_t,~~0\le t\le T,~~\dbP\mbox{-a.s.}
\eea
\end{defn}
\noindent The above $\pa_t u$, $\pa_\o u$ and $\pa^2_{\o\o} u$, if they exist, are unique. Consequently, we call them the time derivative, the first order and second order space derivatives of $u$, respectively.

\begin{defn}
\label{defn-classical}  We say $u\in C^{1,2}(\L)$ is a classical solution (resp. supersolution, subsolution) of PPDE \reff{PPDE} if $\cL u (t,\o) =$ (resp. $\ge, \le$) $0$, for all $(t,\o) \in [0, T)\times \O$.  
\end{defn}
\subsubsection{The shifted spaces}
 Fix  $0\le t\le T$.

-  Let $\O^t:= \big\{\o\in C([t,T], \dbR^d): \o_t ={\bf 0}\big\}$ be the shifted canonical space; $B^{t}$ the shifted canonical process on
$\O^t$;   $\dbF^{t}$ the shifted filtration generated by $B^{t}$, $\dbP^t_0$ the Wiener measure on $\O^t$, and $\L^t :=  [t,T]\times \O^t$.

- For $s\in [t, T]$, define $\|\cd\|_s$ on $\O^t$ and $\dbf$ on $\L^t$ in the spirit of (\ref{rho}),  and the sets $\dbL^0(\L^t)$ etc. in an obvious way. 

- For  $s\in [0, t]$, $\o\in \O^s$ and $\o'\in \O^t$, define the concatenation path $\o\otimes_{t} \o'\in \O^s$ by:
\beaa
(\o\otimes_t \o') (r) := \o_r\1_{[s,t)}(r) + (\o_{t} + \o'_r)\1_{[t, T]}(r),
&\mbox{for all}&
r\in [s,T].
\eeaa

- Let $s\in[0,T)$,  $\xi \in \dbL^0(\cF^s_T)$, and $X\in \dbL^0(\L^s)$. For $(t, \o) \in \L^s$, define $\xi^{t,\o} \in \dbL^0(\cF^t_T)$ and $X^{t,\o}\in \dbL^0(\L^t)$ by:
\beaa
\xi^{t, \o}(\o') :=\xi(\o\otimes_t \o'), \q X^{t, \o}(\o') := X(\o\otimes_t \o'),
&\mbox{for all}&
\o'\in\O^t.
\eeaa
Moreover, for a random time $\t$, we shall use the notation $\xi^{\t,\o} := \xi^{\t(\o),\o}$.

- Define $\cT^t$, $\cH^t$, $\cP^t_L$, $\cP_\infty^t$, and $C^{1,2}(\L^t)$ etc. in an obvious manner. 

It is clear that $u^{t,\o} \in C^0(\L^t)$ for any  $u\in C^0(\L)$ and $(t,\o) \in \L$. Similar property holds for other spaces introduced above. 
Moreover, for any $\t\in \cT$ (resp. $\ch\in\cH$) and any $(t,\o)\in\L$ such that $t<\t(\o)$ (resp. $t<\ch(\o)$), it is clear that $\t^{t,\o}\in \cT^t$ (resp. $\ch^{t,\o}\in\cH^t$).

\subsubsection{Viscosity solutions of PPDEs}

We first introduce the spaces for viscosity solutions. 
\begin{defn}
\label{defn-spaceC0} Let  $u \in \dbL^0 (\L)$.
\\
(i)  We say $u$ is right continuous in $(t,\o)$ under $\dbf$ if: for any $(t,\o) \in \L$ and any $\e>0$, there exists $\d>0$ such that, for any $(s, \tilde\o)\in \L^t$ satisfying $\dbf((s,  \tilde\o), (t, {\bf 0})) \le \d$, we have $|u^{t,\o}(s,\tilde \o) - u(t,\o)|\le \e$. 
\\
(ii) We say $u\in \Usub$ if $u$ is bounded from above, right continuous in $(t,\o)$ under $\dbf$, and  there exists a modulus of continuity function $\rho$ such that for any $(t,\o), (t',\o')\in\L$:
 \bea\label{USC}
 u(t,\o) - u(t', \o') 
 \le 
 \rho\Big(\dbf\big((t,\o), (t',\o')\big)\Big)
 ~\mbox{whenever}~t\le t'.
 \eea 
(iii) We say $u\in \Usup$ if $-u \in\Usub$.
\end{defn}
It is clear that $\Usub ~\!\cap~\! \Usup = UC_b(\L)$. We also recall from \cite{ETZ0} Remark 3.2 that Condition \reff{USC} implies that $u$ has left-limits and positive jumps.

We next introduce the  nonlinear expectations. Denote by $\dbL^1(\cF^t_T,\cP^t_L)$ the set of $\xi\in \dbL^0(\cF^t_T)$ with $\sup_{\dbP\in\cP^t_L}\dbE^{\dbP}[|\xi|]<\infty$, and define, for $\xi\in \dbL^1(\cF^t_T,\cP^t_L)$,
 \beaa
 \overline{\cE}^L_t[\xi]
 =
 \sup_{\dbP\in\cP^t_L}\dbE^{\dbP}[\xi]
 ~~\mbox{and}~~
 \underline{\cE}^L_t[\xi]
 =
 \inf_{\dbP\in\cP^t_L}\dbE^{\dbP}[\xi]
 =
 -\overline{\cE}^L_t[-\xi].
 \eeaa

We now define viscosity solutions. 
For any  $u\in \dbL^0(\L)$, $(t,\o) \in [0, T)\times \O$, and  $L>0$, let
 \begin{equation}\label{cA}
 \ba{lll}
\displaystyle \!\!\underline\cA^{\!L}\!u(t,\o) 
\! :=\!
 \Big\{\f\in C^{1,2}(\L^{\!t})\!\!: \exists \ch\in \cH^t~\mbox{s.t.}~
       (\f-u^{t,\o})_t 
       = 0 =
      \inf_{\t\in \cT^t} \underline\cE^L_t\big[(\f-u^{t,\o})_{\t\wedge\ch}
                       \big]
   
 \Big\},
 \\
\displaystyle \!\! \overline\cA^{\!L}\!u(t,\o) 
 \!:=\! 
 \Big\{\f \in C^{1,2}(\L^{\!t})\!\!: \exists \ch\in \cH^t~\mbox{s.t.}~
      (\f-u^{t,\o})_t 
      =0=
     \sup_{\t\in \cT^t}    \overline\cE^L_t\big[(\f-u^{t,\o})_{\t\wedge\ch}
                       \big] 
 \Big\}.
 \ea
 \end{equation}

\begin{defn}
\label{defn-viscosity}
{\rm (i)} Let $L>0$. We say $u\in\Usub$ (resp. $\Usup$) is a viscosity $L$-subsolution (resp. $L$-supersolution) of PPDE (\ref{PPDE})  if,  for any $(t,\o)\in [0, T)\times \O$ and any $\f \in \underline\cA^{L}u(t,\o)$ (resp. $\f \in \overline\cA^{L}u(t,\o)$):
 \bea
 \label{cLto}
\cL^{t,\o}\f (t,{\bf 0}) :=  \Big[-\pa_t \f -G^{t,\o}(., \f,\pa_\o \f,\pa^2_{\o\o}\f) 
 \Big](t,{\bf 0}) 
 &\le  ~~(\mbox{resp.} \ge)&  0.
 \eea

\noindent {\rm (ii)} We say $u\in\Usub$ (resp. $\Usup$)  is a viscosity subsolution (resp. supersolution) of PPDE (\ref{PPDE}) if  $u$ is viscosity $L$-subsolution (resp. $L$-supersolution) of PPDE (\ref{PPDE}) for some $L>0$.

\noindent {\rm (iii)} We say $u\in{\rm UC}_b(\L)$ is a viscosity solution of PPDE (\ref{PPDE}) if it is both a viscosity subsolution and a viscosity supersolution.
\end{defn}
As pointed out in \cite{ETZ1} Remark 3.11 (i), without loss of generality in \reff{cA}  we may always set $\ch = \ch^t_\e$ for some small $\e>0$:
\bea
\label{chet}
\ch^t_\e := \inf\{s> t: |B^t_s|\ge \e\} \wedge (t+\e).
\eea

\subsection{Monotone schemes for (standard) PDEs}
In this subsection we introduce the main result of  Barles and Souganidis  \cite{BS}. We shall follow the presentation in Guo, Zhang and Zhuo \cite{GZZ}.
We first recall the definition of viscosity solutions for PDE \reff{PDE}:  an upper (resp. lower) semicontinuous function $v$  is called a viscosity subsolution (resp. viscosity supersolution) of PDE \reff{PDE} if $\dbL \f (t,x) \le (\mbox{resp.} \ge )~ 0$,
 for any $(t,x) \in [0, T)\times \dbR^d$ and any smooth function $\f$ satisfying:
\bea
\label{PDE-test}
[u-\f](t,x) = 0 \ge (\mbox{resp.} \le) [u-\f](s,y),\q \mbox{for all}~(s,y) \in [0, T]\times \dbR^d.
\eea
For the viscosity theory of PDEs, we refer to the classical references  \cite{CIL, FS, YZ}. 

We shall adopt the following standard assumptions:

\begin{assum}
\label{assum-PDE}
(i)  $F(\cd, 0, {\bf 0}, {\bf 0})$ and $f$ are bounded.

(ii) $F$ is continuous in $t$, uniformly Lipschitz continuous in $(x, y, z, \g)$, and $f$ is uniformly Lipschitz continuous in $x$.

(iii) 
PDE \reff{PDE} is parabolic, that is, $F$ is nondecreasing in $\g$.

(iv)  Comparison principle for PDE \reff{PDE}  holds in the class of bounded viscosity solutions. That is, if $v_1$ and $v_2$ are bounded viscosity subsolution and viscosity supersolution of PDE \reff{PDE}, respectively, and $v_1(T,\cd) \le f\le v_2(T,\cd)$, then $v_1 \le v_2$ on $[0, T]\times \dbR^d$.
\end{assum}
For any $t\in [0, T)$ and $h \in (0, T-t)$, let $\dbT^{t,x}_h$ be an operator on the set of measurable functions $\f :  \dbR^d \to \dbR$. For $n\ge 1$, denote $h:= {T\over n}$, $t_i := i h$, $i=0,1,\cds, n$, and define:
\bea
\label{PDE-vh}
v^h(t_n, x) := f(x),~~ v^h(t, x) := \dbT^{t,x}_{t_i - t}[v^h(t_{i},\cd)],~ t \in [t_{i-1}, t_i),~ i=n,\cds, 1.
\eea
The following convergence result is reported in \cite{GZZ} Theorem 2.2, which is based on \cite{BS} and is due to Fahim, Touzi and Warin \cite{FTW} Theorem 3.6.
\begin{thm}
\label{thm-PDE}
Let Assumption \ref{assum-PDE} hold. Assume  $\dbT^{t,x}_h$ satisfies the following  conditions:

(i) Consistency:   for any $(t,x)\in [0, T)\times \dbR^d$ and any $\f\in C^{1,2}([0, T)\times \dbR^d)$,
\beaa
\lim_{(t', x', h, c) \to (t,x,0,0)}{[c+\f](t',x')- \dbT^{t',x'}_h\big[[c+\f](t'+h, \cd)\big]\over h} = \dbL\f(t,x).
\eeaa

(ii) Monotonicity: $\dbT^{t,x}_h [\f] \le \dbT^{t,x}_h[\psi]$ whenever $\f \le \psi$.

(iii) Stability:  $v^h$ is bounded uniformly in $h$ whenever $f$ is bounded.

(iv) Boundary condition: $\lim_{(t',x',h)\to (T,x,0)}v^h(t',x')= f(x)$  for any $x\in\dbR^d$.

\noindent Then PDE (\ref{PDE}) with terminal condition $v(T,\cd) = f$ has a unique bounded viscosity solution $v$, and $v^h$ converges to  $v$ locally uniformly as $h\to 0$.
 \end{thm}

\section{Monotone scheme for PPDEs}
\label{sect-mon}
\setcounter{equation}{0}
Our goal of this section is to extend Theorem \ref{thm-PDE} to PPDE \reff{PPDE}. Similar to Assumption \ref{assum-PDE}, we assume
 \begin{assum}
\label{assum-PPDE}
(i)   $G(\cd, 0, {\bf 0}, {\bf 0})$ and $g$ are bounded.

(ii) $G$ is  continuous in $(t,\o)$, uniformly Lipschitz continuous in $(y, z, \g)$, and $g$ is uniformly continuous in $\o$. Denote by $L_0$ the Lipschitz constant of $G$ in $(z,\g)$. 

(iii) 
PDE \reff{PDE} is parabolic, that is, $G$ is nondecreasing in $\g$.

(iv)  Comparison principle for PPDE \reff{PPDE}  holds in the class of bounded viscosity solutions. That is, if $u_1$ and $u_2$ are bounded viscosity subsolution and viscosity supersolution of PPDE \reff{PPDE}, respectively, and $u_1(T,\cd) \le g\le u_2(T,\cd)$, then $u_1 \le u_2$ on $\L$.
\end{assum}
For the comparison principle in (iv) above, we refer to \cite{ETZ2} for some sufficient conditions. 

Now for any $(t,\o)\in [0, T)\times\O$ and $h \in (0, T-t)$, let $\dbT^{t,\o}_h$ be an operator on $\dbL^0(\cF^t_{t+h})$.  For $n\ge 1$, denote $h:= {T\over n}$, $t_i := i h$, $i=0,1,\cds, n$, and define:
\bea
\label{PPDE-uh}
u^h(t_n, \o) := g(\o),~~ u^h(t, \o) := \dbT^{t,\o}_{t_i - t}\big[u^h(t_{i},\cd)\big],~ t \in [t_{i-1}, t_i),~ i=n,\cds, 1.
\eea
where we abuse the notation that: 
\beaa
\dbT^{t,\o}_h[\f] := \dbT^{t,\o}_h[\f^{t,\o}],\q \mbox{for} ~\f\in \dbL^0(\cF_{t+h}).
\eeaa
The following main result is analogous to Theorem \ref{thm-PDE}. 
\begin{thm}
\label{thm-PPDE}
Let Assumption \ref{assum-PPDE} hold. 
Assume $\dbT^{t,\o}_h$ satisfies the following  conditions:

(i) Consistency:   for any $(t,\o)\in [0, T)\times \O$ and $\f\in C^{1,2}(\L^t)$,
\bea
\label{consistency}
\lim_{(t', \o', h, c) \to (t,{\bf 0},0,0)}{[c+\f](t',\o')- \dbT^{t', \o\otimes_t \o'}_h\big[[c+\f](t'+h, \cd)\big]\over h} = \cL^{t,\o}\f(t,{\bf 0}).
\eea
where $(t', \o')\in \L^t$, $h\in (0, T-t)$, $c\in \dbR$, and $\cL^{t,\o}\f$ is defined in \reff{cLto}.

(ii) Monotonicity: for some constant $L\ge L_0$ and any $\f, \psi \in UC_b(\cF^t_{t+h})$, 
\bea
\label{mon}
\ol \cE^{L}_t[ \f -\psi] \le 0\q\mbox{implies}\q \dbT^{t,\o}_h[\f] \le \dbT^{t,\o}_h[\psi].
\eea

(iii) Stability:  $u^h$ is uniformly bounded and uniformly continuous in $\o$, uniformly on $h$. Moreover, there exists a modulus of continuity function $\rho$, independent of $h$,  such that
\bea
\label{treg}
|u^h(t,\o) - u^h(t', \o_{\cd\wedge t})|\le \rho\Big((t'-t)\vee h\Big),\q\mbox{for any $t<t'$ and any $\o\in \O$}.
\eea



\noindent Then PPDE (\ref{PPDE}) with terminal condition $u(T,\cd) = g$ has a unique bounded $L$-viscosity solution $u$, and $u_h$ converges to  $u$ locally uniformly as $h\to 0$.
 \end{thm}
 
 \begin{rem}
 \label{rem-PPDE}
{\rm  The conditions in Theorem \ref{thm-PPDE} reflect the features of our definition of viscosity solution for PPDEs.

(i) For the consistency condition \reff{consistency}, we require the convergence only for $t'\ge t$.

(ii) The monotonicity condition in Theorem \ref{thm-PDE} (ii) is due to the maximum condition \reff{PDE-test} in the definition of viscosity solutions for PDEs. In our path dependent case,  the monotonicity condition  \reff{mon} is modified in a way to adapt to \reff{cA}.

(iii) Due to the uniform continuity required in the definition of viscosity solutions, the stability condition in Theorem \ref{thm-PPDE} (iii) is somewhat strong. Note that this condition obviously implies the counterparts of the Stability and Boundary conditions in Theorem \ref{thm-PDE}. 
\qed}
\end{rem}
To prove the theorem, we need a technical lemma.

\begin {lem}
\label{lem-ch}
Let $L>0$, $\ch\in \cH$, $\t\in \cT$, $\t \le \ch$, and $X\in \Usub$ with modulus of continuity function $\rho$ in \reff{USC}. Assume
\bea
\label{chcondition}
\ol{\cE}_0^L[{X}_{\t}]-\ol{\cE}_0^L\big[X_\ch\big]  \ge c>0
 \eea
Then there exist constants $\d_0 = \d_0(c,L, d,\rho)>0$,  $C = C(L, d)>0$, and  $\o^*\in \O$ such that  
 \bea
 \label{chest}
t_*:= \t(\o^*)<\ch(\o^*)\q\mbox{and}\q\sup_{\dbP\in \cP^{t_*}_L}\dbP \big[\ch^{t_*, \o^*} -t_* \leq \d \big ]\leq C \d^2~~\mbox{for all}~\d\le \d_0.
  \eea
\end{lem}

\proof Let $\ch$ correspond to $O$ and $t_0$ in \reff{ch}. We first claim  there exist $\d_0=\d_0(c,L, d,\rho)$ and $\o^*$ such that
\bea
\label{oh}
t_*:= \t(\o^*)<t_0-\d_0\q\mbox{and}\q d(\o^*_{t_*}, O^c) \ge \d_0^{1\over 6}.
\eea
In particular, this implies that $t_*<\ch(\o^*)$.  Then, for any $\dbP\in \cP^{t_*}_L$ and $\d\le \d_0$, 
\beaa
&&\dbP \big[\ch^{t_*, \o^*} -t_* \leq \d \big ] =\dbP \Big(\ch^{t_*, \o^*} -t_* \leq \d,~ \o^*_{t_*}+B^{t_*}_{ \ch^{t_*, \o^*}} \in O^c\Big)\\
&\le& \dbP \Big( \sup_{t_* \le s \le t_*+\d} |B^{t_*}_s| \ge d(\o^*_{t_*}, O^c)\Big) \le  \dbP \Big( \sup_{t_* \le s \le t_*+\d} |B^{t_*}_s| \ge \d_0^{1\over 6}\Big)\\
&\le& \d_0^{-1}\dbE^\dbP\Big[ \sup_{t_* \le s \le t_*+\d} |B^{t_*}_s|^6\Big] \le C\d^2,
\eeaa
proving \reff{chest}.

We now prove \reff{oh} by contradiction. Assume \reff{oh} is not true, then
\bea
\label{oh2}
\t \ge t_0-{\d_0} \q\mbox{or}\q  d(B_\t, O^c) < {\d_0}^{1\over 6},\q \forall \o\in\O.
\eea
By definition of $\ol\cE_0^L,$ there exists $\dbP\in \cP_L^0$ such that
\bea\label{specificP}
\ol\cE_0^L[X_\tau] \leq \dbE^\dbP [X_\tau]+{c\over 2}.
\eea
Note that  $B_\t(\o) \in O$ whenever $\t(\o)<\ch(\o)$. Recall \reff{cPL} and let $\eta(\o)$ denote the unit vector pointing from $B_{\t}(\o)$ to $O^c$. Set $\eta(\o)$ be a fixed unit vector when $\t(\o) = \ch(\o)$. Then $\eta \in \cF_\t$. Construct $\hat{\dbP}\in \cP_L^0$ as follows:
\beaa
\a^{\hat\dbP}_t := \a^\dbP_t \1_{[0, \t)}(t) +L \eta \1_{[\t, t_0)},\q \b^{\hat\dbP}_t := \b^\dbP_t \1_{[0, \t)}(t).
\eeaa
That is, $\hat \dbP = \dbP$ on $\cF_\t$ and $dB^{\t(\o)}_t = L\eta(\o) dt$, $t\ge \t$, $\hat\dbP^{\t,\o}$-a.s., where   $\hat\dbP^{\t,\o}$ is the regular conditional probability distribution of $\dbP$. 
 Then, one can easily see that 
\beaa
|B^{\t(\o)}_t| = L[t-\t(\o)],\q \ch^{\t,\o} - \t(\o) = {d(B_\t(\o), O^c)\over L} \wedge [t_0-\tau(\o)],~~\hat{\dbP}^{\tau, \o}\mbox{-a.s. for all}~\o.
\eeaa
This, together with \reff{oh2}, implies
\bea
\label{hatPdist}
\dbf\big((\t, \o), (\ch, \o)\big) = \ch-\t + \sup_{\t\le t\le \ch} |B^\t_t| \le C[\ch-\t] \le C[{{\d_0}^{1\over 6}\over L} + {\d_0}]\le C{\d_0}^{1\over 6}, \q\hat{\dbP}^{\tau, \o}\mbox{-a.s.}
\eea
 Then, by \reff{specificP}, \reff{USC}, and \reff{hatPdist},
\beaa
 && \ol{\cE}_0^L[{X}_{\t}]-\ol{\cE}_0^L\big[X_\ch\big] \le \dbE^\dbP[X_\t] - \dbE^{\hat\dbP}[X_\ch] +{c\over 2} =  \dbE^{\hat\dbP}[X_\t-X_\ch] +{c\over 2}\\
  &\le& \dbE^{\hat\dbP}\Big[\rho\Big(\dbf\big((\t, \o), (\ch, \o)\big)\Big)\Big] +{c\over 2} \le \rho\big(C\d_0^{1\over 6}\big)+{c\over 2}.
  \eeaa
This contradicts with \reff{chcondition} when ${\d_0}$ is small enough, and thus \reff{oh} holds true.
\qed

\bs

\noindent{\bf Proof of Theorem \ref{thm-PPDE}.} By the stability, $u^h$ is bounded. Define
\bea
\label{ulu}
\ulu\tow:=\liminf_{h\rightarrow 0} u^h\tow,\q \olu\tow:=\limsup_{h\rightarrow 0} u^h\tow.
\eea
Clearly $\ul u (T,\o) = g(\o) = \ol u(T,\o)$, $\ul u\le \ol u$, and $\ul u, \ol u$ are bounded and uniformly continuous. We shall show that $\ul u$ (resp. $\ol u$) is a viscosity $L$-supersolution (resp. $L$-subsolution) of PPDE \reff{PPDE}. Then by the comparison principle we see that $\ol u \le \ul u$ and thus $u:=\ol u = \ul u$ is the unique viscosity solution of PPDE \reff{PPDE}.   The convergence of $u^h$ is obvious now, which, together with the uniform regularity of $u^h$ and $u$, implies further the locally uniform convergence.

Without loss of generality, we shall only prove  by contradiction that $\ul u$  satisfies the viscosity $L$-supersolution  property  at $(0, {\bf 0})$. Assume not, then there exists $\f^0 \in \ol \cA^{L} \ul u(0,{\bf 0})$ with corresponding $\ch \in \cH$ such that
$
-c_0:=\cL \f^0(0,\zo)<0.
$
Denote 
\bea
\label{f0}
\f(t,\o) := \f^0(t,\o) - {c_0\over 2}t.
\eea
Then
\bea
\label{cLf0}
\cL \f(0,\zo) = -{c_0\over 2}<0.
\eea

Denote $X^0:=\f-\ulu$, $X^h := \f-u^h$, and $\ol \cE:= \ol \cE^{L}_0$, $\ul \cE:= \ul \cE^{L}_0$. Recall \reff{chet} and denote $\ch_\e := \ch^0_\e \wedge \e^5$, $c_\e := {1\over 3}c_0\e^5$. Note that $\ch_\e  \le \ch$  for $\e$ small enough, and by \cite{ETZ1} (2.8),
\bea
\label{cheeta}
\sup_{\dbP\in \cP_{L}} \dbP(\ch_\e \neq \e^5) = \sup_{\dbP\in \cP_{L}}  \dbP(\ch_\e^0 < \e^5) \le CL^4 \e^{-4} \e^{10} \le C\e c_\e .
\eea
Then
\beaa
\ol\cE[\e^5 - \ch_\e] \le \ol\cE\big[\e^5 \1_{\{ \ch_\e\neq \e^5\}}\big] \le C\e c_\e.
\eeaa
  Thus, for $\e$ small, it follows from $\f^0\in \ol\cA^{L} \ul u(0,{\bf 0})$ that
\bea
\label{ce}
X^0_0 - \ol \cE[X^0_{\ch_\e}] &=& [\f^0-\ul u]_0 -  \ol \cE\Big[(\f^0-\ul u)_{\ch_\e} - {c_0\over 2} {\ch_\e}\Big] \nonumber\\
&\ge&  \ol \cE\Big[(\f^0-\ul u)_{\ch_\e} \Big] - \ol \cE\Big[(\f^0-\ul u)_{\ch_\e} - {c_0\over 2} {\ch_\e}\Big] \\
&\ge&  \ul \cE\Big[{c_0\over 2} {\ch_\e}\Big]  = {c_0\e^5\over 2} - {c_0\over 2} \ol\cE[\e^5-{\ch_\e}] \ge {3c_\e\over 2} - C\e c_\e \ge c_\e >0.\nonumber
\eea
Let $h_k\downarrow 0$ be a sequence such that
\bea
\label{hk}
\lim_{k\to\infty} u^{h_k}_0 = \ul u_0,
\eea
and simplify the notations: $u^k := u^{h_k}$, $X^k := X^{h_k}$.  Then \reff{ce} leads to
\beaa
c_\e \le [\f_0-\liminf_{h\to 0} u^h_0] -  \ol \cE\Big[\f_{\ch_\e} - \liminf_{h\to 0} u^h_{\ch_\e} \Big]  \le [\f_0-\lim_{k\to \infty} u^{k}_0] -  \ol \cE\Big[\f_{\ch_\e} - \liminf_{k\to\infty} u^{k}_{\ch_\e} \Big].
\eeaa
Note that $X^k$ is uniformly bounded. Then by \reff{cheeta} we have
\beaa
\ol\cE\Big[|X^k_{\ch_\e}-X^k_{\e^5}|\Big] \le C\e c_\e.
\eeaa
Since $u^h$ is uniformly continuous, applying the monotone convergence theorem under nonlinear expectation $\ol\cE$, see e.g. \cite{ETZ0} Proposition 2.5, we have
\beaa
c_\e  & \le& \lim_{k\to \infty} [\f_0- u^{k}_0] -  \ol \cE\Big[ \limsup_{k\to\infty} [\f_{\ch_\e} - u^{k}_{\ch_\e}] \Big] \\
&\le& \lim_{k\to \infty} X^k_0 -  \ol \cE\Big[ \limsup_{k\to\infty} X^k_{\e^5} \Big] + C\e c_\e =  \lim_{k\to \infty} X^k_0 -  \ol \cE\Big[ \lim_{m\to\infty}\sup_{k\ge m} X^k_{\e^5} \Big] + C\e c_\e\\
&=& \lim_{k\to \infty} X^k_0 -   \lim_{m\to\infty}\ol \cE\Big[\sup_{k\ge m} X^{k}_{\e^5}\Big]+ C\e c_\e \le  \lim_{k\to \infty} X^k_0 -   \limsup_{k\to\infty}\ol \cE\Big[X^{k}_{\e^5}\Big]+ C\e c_\e \\
&\le& \lim_{k\to \infty} X^{k}_0 -  \limsup_{k\to\infty} \ol \cE\Big[X^{k}_{\ch_\e} \Big]+ C\e c_\e=  \liminf_{k\to\infty}\Big[ X^{k}_0-\ol \cE\big[X^{k}_{\ch_\e} \big]\Big] + C\e c_\e.
\eeaa
Then, for all $\e$ small enough and  $k$ large enough,
\bea
\label{Xk0}
X^{k}_0-\ol \cE\big[X^{k}_{\ch_\e} \big] \ge {c_\e\over 2}.
\eea

Now for each $k$, define
\beaa
Y^k_t(\o) := \sup_{\t\in \cT^t}\ol \cE^{L}_t\big[(X^{k})^{t,\o}_{\t\wedge \ch_\e^{t,\o}}\big],~ t\le {\ch_\e}(\o),\q\mbox{and}\q \t_k := \inf\{t\ge 0: Y^k_t = X^{k}_t\}.
\eeaa 
We remark that here $Y^k, \t_k$ depend on $\e$ as well, but we omit the superscript $^\e$ for notational simplicity. 
Applying \cite{ETZ0} Theorem 3.6, we know $\t_k\le {\ch_\e}$ is an optimal stopping time for $Y^k_0$ and thus 
\beaa
0<{c_\e\over 2} &\le& X^{k}_0-\ol \cE\big[X^{k}_{\ch_\e} \big] \le Y^k_0 -\ol \cE\big[X^{k}_{\ch_\e} \big] = \ol \cE\big[X^{k}_{\t_k} \big]-\ol \cE\big[X^{k}_{\ch_\e} \big]
\eeaa
By Lemma \ref{lem-ch}, for $k$ large enough so that $h_k\le \d_0({c_\e\over 2}, L, d, \rho)$, there exists $\o^k\in \O$ such that 
\bea
\label{tk*}
t^k_* := \t_k(\o^k) < {\ch_\e}(\o^k) \q\mbox{and\q $\sup_{\dbP\in \cP^{t_*^k}_{L}}$}~ \dbP\Big(\ch_\e^{k}- t_*^k \le \d\Big) \le C\d^2~~\mbox{for all}~\d\le h_k,
\eea
where $\ch_\e^k := \ch_\e^{t_*^k, \o^k}$.
Let $\{t^k_i, i=0,\cds, n_k\}$ denote the time partition corresponding to $h_k$, and assume $t^k_{i-1} \le t^k_* < t_i^k$. Note that 
\beaa
X^k_{t^k_*} (\o^k) = Y^k_{t^k_*} (\o^k)  \ge \overline\cE^{L}_{t^k_*}\Big[(X^k)^{t^k_*,\o^k}_{\t\wedge \ch_\e^k}\Big],\q\forall \t \in \cT^{t^k_*}.
\eeaa
Set $\d_k := t^k_i - t^k_* \le h_k$ and $\t:= t^k_i$. Combine the above inequality and \reff{tk*} we have
\beaa
[\f - u^k](t^k_*,\o^k) \ge \overline\cE^{L}_{t^k_*}\Big[(\f-u^k)^{t^k_*,\o^k}_{t^k_i\wedge \ch_\e^k}\Big] \ge \overline\cE^{L}_{t^k_*}\Big[(\f-u^k)^{t^k_*,\o^k}_{t^k_i}\Big] - C \d_k^2.
\eeaa
This implies
\beaa
 \overline\cE^{L}_{t^k_*}\Big[\Big(\f^{t^k_*,\o^k}_{t^k_i} -  [\f - u^k](t^k_*,\o^k) - C \d_k^2\Big) - (u^k)^{t^k_*,\o^k}_{t^k_i}\Big] \le 0.
 \eeaa
 By the monotonicity condition \reff{mon} we have
 \bea
 \label{uktk}
 u^k(t^k_*,\o^k) = \dbT^{t^k_*,\o^k}_{\d_k}[u^k_{t^k_i}]  \le  \dbT^{t^k_*,\o^k}_{\d_k}\Big[\f_{t^k_i} -  [\f - u^k](t^k_*,\o^k) - C \d_k^2\Big].
 \eea
 
We next use the consistency condition \reff{consistency}. For $(t,\o) = (0, {\bf 0})$, set
\beaa
t':= t^k_*,\q \o' := \o^k,\q h:= \d_k,\q c := -  [\f - u^k](t^k_*,\o^k) - C \d_k^2.
\eeaa
By first sending $k\to \infty$ and then $\e\to 0$, we see that
\beaa
\dbf((t^k_*, \o^k), (0,{\bf 0})) \le \ch_\e + \sup_{0\le t\le \ch_\e} |\o^k_t|\le 2\e \to 0,\q h \le h_k\to 0,
\eeaa
which, together with \reff{f0}, \reff{hk}, and the uniform continuity of $\f$ and $u^k$, implies
\beaa
|c| \le \Big| [\f - u^k](t^k_*,\o^k) -  [\f - u^k](0, {\bf 0})\Big| + | u^k_0 - \ul u_0| + C\d_k^2 \to 0.
\eeaa
Then, by the consistency condition \reff{consistency} we obtain from \reff{uktk} that
\beaa
0 &\le& {u^k(t^k_*,\o^k)  - \dbT^{t^k_*,\o^k}_{\d_k}\Big[\f_{t^k_i} -  [\f - u^k](t^k_*,\o^k) - C \d_k^2\Big]\over \d_k} \\
&=&{ [c + \f](t^k_*, \o^k) - \dbT^{t^k_*,\o^k}_{\d_k}\Big[ [c+\f]_{t^k_i}\Big]\over \d_k} + C\d_k\to  \cL \f(0, {\bf 0}).
\eeaa
This contradicts with \reff{cLf0}. 
\qed

\section{An illustrative monotone scheme}
\label{sect-scheme}
\setcounter{equation}{0}

We first remark that the monotonicity condition \reff{mon} is solely due to our definition of viscosity solution of PPDEs. It is sufficient but not necessary for the convergence of the scheme. In Markovian case, the PPDE \reff{PPDE} is reduced back to PDE \reff{PDE}. The schemes proposed in \cite{FTW} and \cite{STZ} satisfy the traditional monotonicity condition in Theorem \ref{thm-PDE}, but violates our new monotonicity condition \reff{mon}. However, as proved in \cite{FTW, STZ}, we know those schemes do converge.

The goal of this section is to propose a scheme which satisfies all the conditions in Theorem \ref{thm-PPDE} and thus converges. However, to ensure the monotonicity condition \reff{mon}, we will need certain conditions which are purely technical. Monotone schemes for general parabolic PPDEs is a challenging problem and we shall leave it for future research. We also remark that  efficient implementation of such schemes, especially in high dimensions, is also a very challenging problem and will also be left for future research.

Our scheme will involve some parameters:
\bea
\label{parameter}
\mu_i >0, ~ \si_i>0,~ i=1,\cds, d.
\eea
Let $e_i\in \dbR^d$ be the vector whose $i$-th component is $1$ and all other components are $0$, and $e_{ij}\in \dbR^{d\times d}$ be the matrix whose $(i,j)$-th component is $1$ and all other components are $0$.  
Given $(t,\o)\in [0, T)\times \O$, recall \reff{cPL} and introduce the following probability measures on $\O^t$: for $i, j=1,\cds,d$, 
\bea
\label{dbPi}
\left.\ba{c}
\dbP^0:\q  \a^{\dbP^0} = {\bf 0},  \b^{\dbP^0}={\bf 0};\qq \dbP^i: \q \a^{\dbP^i} = \mu_i e_i,  \b^{\dbP^i} = {\bf 0};\\
 \dbP^{ii}:\q \a^{\dbP^{ii}} = {\bf 0},  \b^{\dbP^{ii}} = \si_i e_{ii};\qq   \dbP^{ij}:\q \a^{\dbP^{ij}} = {\bf 0},  \b^{\dbP^{ij}} = \si_je_{ij} + \si_i e_{ji}, ~i\neq j.
 \ea\right.
\eea
Now for $h\in (0, T-t)$ and $\f\in \dbL^0(\cF^t_{t+h})$, define
\bea
\label{dbTf}
\dbT^{t,\o}_h[\f] := \cD^{(0)} \f+ h G(t, \o,\cD^{(0)} \f, \cD^{(1)} \f, \cD^{(2)} \f),
\eea
where $\cD^{(0)} \f$, $ \cD^{(1)} \f$, $\cD^{(2)} \f$ take values in $\dbR$, $\dbR^d$, $\dbS^d$, respectively, with each component defined by
\bea
\label{cDf}
\left.\ba{c}
\displaystyle \cD^{(0)} \f := \dbE^{\dbP^0}[\f],\q \cD^{(1)}_i \f := {\dbE^{\dbP^i}[\f] - \dbE^{\dbP^0}[\f]\over \mu_i h}, \q\cD^{(2)}_{i,i} \f := {\dbE^{\dbP^{ii}}[\f]-  \dbE^{\dbP^0}[\f]\over \si_i^2 h\slash 2},\\
\displaystyle  \cD^{(2)}_{i,j} \f := {\dbE^{\dbP^{ij}}[\f] - \dbE^{\dbP^{ii}}[\f]- \dbE^{\dbP^{jj}}[\f] +\dbE^{\dbP^0}[\f]\over \si_i\si_j h}, \q i\neq j.
\ea\right.
\eea

We now verify the conditions in Theorem \ref{thm-PPDE}. 

\begin{lem} [Consistency]
\label{lem-consistency}
Under Assumption \ref{assum-PPDE}, $\dbT^{t,\o}_h$ satisfies the consistency condition \reff{consistency}.
\end{lem}
\proof Without loss of generality, we assume $(t,\o) = (0, {\bf 0})$. Let $(t',\o', h, c)$ be as in  \reff{consistency}, and for notational simplicity, at below we write $(t', \o')$ as $(t,\o)$. Now for $\f\in C^{1,2}(\L)$, denote $\psi := c + \f^{t,\o}(t+h,\cd) \in \dbL^0(\cF^t_{t+h})$, and $\f|_t^s := \f^{t,\o}(s, B^t) - \f(t,\o)$. Send $(t,\o, h, c)\to (0, {\bf 0}, 0, 0)$, by the functional It\^{o} formula and the smoothness of $\f$, one can easily check that
\beaa
&& \cD^{(0)} \psi = c + \f(t+h,\o_{\cd\wedge t})  \to \f(0,{\bf 0});\\
&&{ \cD^{(0)} \psi- [c+\f](t,\o)\over h} =    {1\over h}\dbE^{\dbP^0}[\f|_t^{t+h}]   = {1\over h}\int_t^{t+h} \pa_t \f(s, \o_{\cd\wedge t})ds \to \pa_t\f(0,{\bf 0});\\
&&\cD^{(1)}_i \psi ={1 \over \mu_i h}\dbE^{\dbP^i}[\f|_t^{t+h}]  -{ 1\over \mu_i h}\dbE^{\dbP^0}[\f|_t^{t+h}] \\
&&\qq = {1\over \mu_i h} \int_t^{t+h} \dbE^{\dbP^i} \Big[\big(\pa_t  + \mu_i \pa_{\o^i}  \big)\f(s,  \o\otimes_t B^t) - \pa_t \f(s, \o_{\cd\wedge t})\Big] ds\to  \pa_{\o^i}\f(0,{\bf 0});\\
&&\cD^{(2)}_{i,i} \psi ={2\over \si_i^2 h}\dbE^{\dbP^{ii}}[\f|_t^{t+h}]  -{2\over \si_i^2 h}\dbE^{\dbP^0}[\f|_t^{t+h}] \\
&&\qq =   {2\over \si_i^2 h}\int_t^{t+h}\dbE^{\dbP^{ii}}\Big[\big(\pa_t  + {\si_i^2\over 2} \pa^2_{\o^i\o^i} \big)\f(s,\o\otimes_t B^t) - \pa_t \f(s,  \o_{\cd\wedge t}) \Big]ds \to  \pa^2_{\o^i\o^i}\f(0,{\bf 0});\\
&&\cD^{(2)}_{i,j} \psi = {1\over \si_i\si_j h} \dbE^{\dbP^{ij}}[\f|_t^{t+h}]  - {1\over \si_i\si_j h} \dbE^{\dbP^{ii}}[\f|_t^{t+h}] - {1\over \si_i\si_j h} \dbE^{\dbP^{jj}}[\f|_t^{t+h}]  + {1\over \si_i\si_j h} \dbE^{\dbP^{0}}[\f|_t^{t+h}]  \\
&&\qq= {1\over \si_i\si_j h}\int_t^{t+h}\dbE^{\dbP^{ij}}\Big[ \Big(\pa_t + {1\over 2} \si_i^2 \pa^2_{\o^i\o^i} + {1\over 2} \si_j^2 \pa^2_{\o^j\o^j} + \si_i\si_j \pa^2_{\o^i\o^j}\Big)\f (s, \o\otimes_t B^t)\Big]ds \\
&&\qq  -{1\over \si_i\si_j h}\int_t^{t+h}\dbE^{\dbP^{ii}}\Big[ \Big(\pa_t + {1\over 2} \si_i^2 \pa^2_{\o^i\o^i} \Big)\f (s, \o\otimes_t B^t)\Big]ds \\
&&\qq - {1\over \si_i\si_j h}\int_t^{t+h}\dbE^{\dbP^{jj}}\Big[ \Big(\pa_t + {1\over 2} \si_j^2 \pa^2_{\o^j\o^j} \Big)\f (s, \o\otimes_t B^t)\Big]ds +  {1\over \si_i\si_j h}\int_t^{t+h} \pa_t\f(s, \o_{\cd\wedge t})ds\\
&&\qq= {1\over \si_i\si_j h}\int_t^{t+h}\dbE^{\dbP^{ij}}\Big[ \Big(\pa_t + {1\over 2} \si_i^2 \pa^2_{\o^i\o^i} + {1\over 2} \si_j^2 \pa^2_{\o^j\o^j} + \si_i\si_j \pa^2_{\o^i\o^j}\Big)\f |_t^s   \Big]ds \\
&&\qq  -{1\over \si_i\si_j h}\int_t^{t+h}\Big(\dbE^{\dbP^{ii}}\Big[ \Big(\pa_t + {1\over 2} \si_i^2 \pa^2_{\o^i\o^i} \Big)\f|_t^s \Big]+\dbE^{\dbP^{jj}}\Big[ \Big(\pa_t + {1\over 2} \si_j^2 \pa^2_{\o^j\o^j} \Big)\f| _t^s \Big]\Big)ds \\
&&\qq+ \pa^2_{\o^i\o^j}\f(t,\o)\to \pa^2_{\o^i\o^j}\f(0, {\bf 0}).
\eeaa
Plug these into \reff{dbTf} and \reff{cDf}, we obtain \reff{consistency} immediately.
\qed

To ensure the monotonicity condition \reff{mon}, we need some additional conditions. 
\begin{assum}
\label{assum-mon}
Assume $G$ is differentiable in $(z, \g)$ and one may choose $\mu_i, \si_i$ so that
\bea
\label{mon1}
\left.\ba{c}
\displaystyle \pa_{z_i} G \ge 0,\q \pa_{\g_{ij}} G \ge 0,\q 2\pa_{\g_{ii}} G\slash \si_i \ge \sum_{j\neq i} [\pa_{\g_{ij}}G + \pa_{\g_{ji}}G]\slash \si_j,\\
\displaystyle \qq \sum_{i=1}^d {\pa_{z_i}G\over \mu_i}  + \sum_{i=1}^d {2\pa_{\g_{ii}} G\over \si_i^2}  - \sum_{i\neq j} {\pa_{\g_{ij}}G\over \si_i\si_j} \le 1 - \e_0\q\mbox{for some}~\e_0\in (0,1).
\ea\right.
\eea
\end{assum}

\begin{rem}
\label{rem-mon}
{\rm (i) The differentiability of $G$ is just for convenience. For notational simplicity, at below we shall assume $G$ is differentiable in $y$ as well.  

(ii) By setting $\si_i$ all equal, a sufficient condition for the third inequality in \reff{mon1} is the following diagonal dominance condition:
\bea
\label{diagonal}
 2\pa_{\g_{ii}} G \ge \sum_{j\neq i} [\pa_{\g_{ij}}G + \pa_{\g_{ji}}G].
 \eea
 
 (iii) Since the derivatives of $G$ are uniformly bounded, thanks to Assumption \ref{assum-PPDE}, then the last inequality in \reff{mon1} always holds true when $\mu_i, \si_i$ are large enough.
\qed}
\end{rem}
\begin{lem} [Monotonicity]
\label{lem-mon}
Under Assumptions \ref{assum-PPDE} and \ref{assum-mon}, $\dbT^{t,\o}_h$ satisfies the monotonicity condition \reff{mon} for $L\ge L_0$ large enough and $h$ small enough.
\end{lem}

\proof  Without loss of generality, we assume $(t,\o) = (0, {\bf 0})$ and denote $\dbT_h := \dbT^{t,\o}_h$.   Assume $L\ge L_0$ is large enough so that the $\dbP^i$ and $\dbP^{ij}$ in \reff{dbPi} are in $\cP_L$. Let $\f_1, \f_2 \in UC_b(\cF_h)$ satisfy 
\bea
\label{cEpsi}
\ol\cE^L[\psi] \le 0,\q\mbox{ where}\q \psi := \f_1-\f_2.
\eea
 Then, recalling \reff{dbTf},
\beaa
\dbT_h\f_1 - \dbT_h \f_2  = \cD^{(0)} \psi + h\Big[\pa_y G \cD^{(0)} \psi + \pa_z G \cd \cD^{(1)} \psi  + \pa_\g G : \cD^{(2)} \psi\Big].
\eeaa
Note that here $\pa_y G$ etc are deterministic.  By \reff{cDf} we have 
\bea
\label{Delta}
\dbT_h\f_1 - \dbT_h \f_2  
&=&a_0 \dbE^{\dbP^0}[\psi] +\sum_{i=1}^d a_i \dbE^{\dbP^i}[\psi]  +  \sum_{i=1}^d a_{ii}\dbE^{\dbP^{ii}}[\psi]  +  \sum_{i\neq j} a_{ij}\dbE^{\dbP^{ij}}[\psi],
\eea
where
\bea
\label{ai}
\left.\ba{c}
\displaystyle a_0:= 1+ h \pa_y G  - \sum_{i=1}^d {\pa_{z_i}G\over \mu_i} - \sum_{i=1}^d {\pa_{\g_{ii}}G \over \si_i^2\slash 2} + \sum_{i\neq j} {\pa_{\g_{ij}}G\over \si_i\si_j}\\
\displaystyle a_i := {\pa_{z_i}G\over \mu_i}, \q a_{ii}:=    {2\pa_{\g_{ii}}G\over \si_i^2} - \sum_{j\neq i} {\pa_{\g_{ij}}G+\pa_{\g_{ji}}G\over \si_i\si_j},\q a_{ij}:= {\pa_{\g_{ij}}G\over \si_i\si_j}.
\ea\right.
\eea
By \reff{mon1} we see that $a_0, a_i, a_{ij} \ge 0$, provided $h$ is small enough. Note that
\beaa
a_0 + \sum_{i=1}^d a_i + \sum_{i,j=1}^d a_{ij} = 1+ h \pa_y G .
\eeaa
Then one may define the following probability measure:
\bea
\label{hatdbP}
\hat\dbP := {1\over 1+h \pa_y G} \Big[a_0 \dbP^0 + \sum_{i=1}^d a_i \dbP^i + \sum_{i, j=1}^d a_{ij} \dbP^{ij}\Big],
\eea
and rewrite \reff{Delta} as
\bea
\label{Delta2}
\dbT_h\f_1 - \dbT_h \f_2  &=& (1+h \pa_y G)\dbE^{\hat\dbP}[\psi]. 
\eea
Since $\dbP^0, \dbP^i, \dbP^{ij} \in \cP_L$,  \reff{cEpsi} implies $\dbE^{\dbP^0}[\psi], \dbE^{\dbP^{i}}[\psi], \dbE^{\dbP^{ij}}[\psi]  \le 0$ and thus $\dbE^{\hat\dbP}[\psi] \le 0$.   This leads to  \reff{mon} immediately.
\qed

\begin{rem}
\label{rem-hatP}
{\rm In general $\hat \dbP$ may not be in $\cP_L$. However, we still have $\ul\cE^L\le \dbE^{\hat\dbP} \le \ol\cE^L$.
\qed}
\end{rem}

We now verify the stability condition.

\begin{lem} [Stability]
\label{lem-stability}
Let Assumptions \ref{assum-PPDE} and \ref{assum-mon} hold, and assume further that $G$ and $g$ are uniformly Lipschitz continuous in $\o$. Then $\dbT^{t,\o}_h$ satisfies the stability condition in Theorem \ref{thm-PPDE} for $L$ large enough and $h$ small enough.
\end{lem}
\proof We assume $L$ and $h$ are chosen so that $\dbT^{t,\o}_h$ satisfies the monotonicity condition \reff{mon}.

(i) We first show that $u^h$ is uniformly bounded. Denote $C_i := C_i^h:= \sup_{\o\in \O} |u^h(t_i, \o)|$, $\f := [u^h(t_{i+1},\cd)]^{t_{i},\o}$, and recall \reff{PPDE-uh}.  By \reff{dbTf}  we have
\beaa
u^h(t_{i},\o) \!=\!   \cD^{(0)} \f+ h G(t_i, \o,\cD^{(0)} \f, \cD^{(1)} \f, \cD^{(2)} \f) - hG(t_i, \o,0, {\bf 0}, {\bf 0}) +hG(t_i, \o,0, {\bf 0}, {\bf 0}).
\eeaa
Following the arguments for \reff{Delta2}, for some $\hat \dbP$ defined in the spirit of \reff{ai}-\reff{hatdbP}, we have 
\bea
\label{rep}
u^h(t_{i},\o) &=& (1+h\pa_y G)  \dbE^{\hat\dbP}[\f]+ h G(t_i, \o,0, {\bf 0}, {\bf 0}).
\eea
Then
\beaa
|u^h(t_{i},\o)| &\le& (1+h\pa_y G)  \Big|\dbE^{\hat\dbP}[\f]\Big| + h|G(t_i, \o,0, {\bf 0}, {\bf 0})|\le (1+Ch) C_{i+1} + Ch.
\eeaa
That is,
\beaa
C_{i} \le [1+Ch] C_{i+1} + Ch.
\eeaa
Note that $C_n = \|g\|_\infty$. Then by the discrete Gronwall inequality we see that $\max_{0\le i\le n} C_i \le C$, where the constant $C$ is independent of $h$. 

Finally, for $t\in (t_{i}, t_{i+1})$, following similar arguments we can easily show that $|u^h(t,\o)|\le [1+Ch] C_{i+1} + Ch\le C$. Therefore, $u^h$ is uniformly bounded.

\ms

(ii) We next show that $u^h$ is uniformly Lipschitz continuous in $\o$.  Let $L_i:= L^i_h$ denote the Lipschitz constant of $u^h(t_i, \cd)$. Given $\o^1, \o^2\in \O$, denote $\psi := [u^h(t_{i+1},\cd)]^{t_{i},\o^1}- [u^h(t_{i+1},\cd)]^{t_{i},\o^2}$, then 
\beaa
|\psi|\le L_{i+1}\|\o^1\otimes_{t_i}B^{t_i}-\o^2\otimes_{t_i}B^{t_i}\|_{t_{i+1}}=L_{i+1}\|\o^1-\o^2\|_{t_{i}}. 
\eeaa
Note that $G$ is uniform Lipschitz continuous in $\o$ with certain Lipschitz constant $L_G$. Then similar to (i) above, we have
\beaa
|u^h(t_{i},\o^1)- u^h(t_{i},\o^2)| &\le& (1+h\pa_y G) \dbE^{\hat\dbP}[|\psi|]  + L_Gh \|\o^1-\o^2\|_{t_{i}} \\
& \le&   L_{i+1}\|\o^1-\o^2\|_{t_{i}} [1+Ch] +  L_Gh \|\o^1-\o^2\|_{t_{i}}.
\eeaa
Then 
\beaa
L_{i} \le L_{i+1} [1+Ch] + L_G h.
\eeaa
 Since $L_n=L_g$ is the Lipschitz constant of $g$ which is independent of $h$, we see that $\max_{0\le i\le n} L_i$ is independent of $h$.  Finally, as in the end of (i) above we see that $u^h(t,\cd)$ is uniformly Lipschitz continuous in $\o$, uniformly in $t$ and $h$.

\ms

(iii) We now prove the following time regularity in two steps: 
\bea
\label{regt}
|u^h(t,\o) - u^h(t', \o_{\cd\wedge t})|\le C\sqrt{t'-t + h},\q \mbox{for all}~0\le t< t'\le T.
\eea

{\it Step 1.} We first assume $t' = T$ and $t = t_i$. For $j=i+1, \cds, n$, in the spirit of \reff{rep}, we may define $\hat \dbP_j$ such that $\hat\dbP_{j+1} = \hat\dbP_{j}$ on $\cF^{t_i}_{t_{j}}$ and  
\beaa
u^h(t_j, \o\otimes_{t_i} B^{t_i}) = [1+ h b_j]\dbE^{\hat \dbP_{j+1}}\Big[u^h(t_{j+1}, \o\otimes_{t_i} B^{t_i}) \Big|\cF^{t_i}_{t_j}\Big] + h c_j,
\eeaa
where $b_j:= \pa_y G(t_j)$ and $c_j:= G(t_j, \o\otimes_{t_i}B^{t_i}, 0,{\bf 0}, {\bf 0})$ are in $\dbL^\infty(\cF^{t_i}_{t_j})$. Denote $\Gamma_i:= 1$, $\Gamma_{j+1}:= \prod_{k=i}^{j}[1+ h b_k]$. By induction we have
\beaa
u^h(t_i, \o) = \dbE^{\hat \dbP_{n}}\Big[ \Gamma_{n}  u^h(t_n, \o\otimes_{t_i} B^{t_i}) +h \sum_{j=i}^{n-1}\Gamma_j c_j \Big] =\dbE^{\hat \dbP_{n}}\Big[ \Gamma_{n}  g(\o\otimes_{t_i} B^{t_i}) +h \sum_{j=i}^{n-1}\Gamma_j c_j \Big]  .
\eeaa
One may easily check that
\beaa
|\Gamma_j |\le C,\q |\Gamma_j - 1| \le C(j-i)h\le C(n-i)h = C(T-t_i).
\eeaa
Thus
\beaa
&&|u^h(t_i,\o) - u^h(t_n,\o_{\cd\wedge t_i})|
= \Big| \dbE^{\hat \dbP_{n}}\Big[ \Gamma_{n}  g(\o\otimes_{t_i} B^{t_i}) +h \sum_{j=i}^{n-1}\Gamma_j c_j  \Big] - g(\o_{\cd\wedge t_i})\Big|\\
&\le&\dbE^{\hat \dbP_{n}}\Big[ |\Gamma_n -1|  |g(\o\otimes_{t_i} B^{t_i})| + |  g(\o\otimes_{t_i} B^{t_i}) - g(\o_{\cd\wedge t_i})|+ C(n-i)h \Big] \\
&\le& C(T-t_i) +  C\dbE^{\hat \dbP_{n}}[\|B^{t_i}\|_T]. 
\eeaa
One can easily show that $\dbE^{\hat \dbP_{n}}[\|B^{t_i}\|_T] \le C\sqrt{T-t_i}$. Then \reff{regt} holds in this case.

{\it Step 2.} We now verify the general case. Assume $t_{i-1} \le t < t_i$ and $t_{j-1} \le t' < t_{j}$, then clearly $i\le j$. Since $u^h(t_{i},\cd)$ and $u^h(t_{j},\cd)$ are Lipschitz continuous in $\o$, by \reff{rep} and following the arguments in Step 1, one can similarly show that
\beaa
|u^h(t,\o) - u^h(t_i, \o_{\cd\wedge t})| &\le& C\sqrt{t_i-t} \le  C\sqrt{h}\\
 |u^h(t',\o_{\cd\wedge t}) - u^h(t_j, \o_{\cd\wedge t})|&\le& C\sqrt{t_j-t'} \le C\sqrt{h},\\
| u^h(t_i, \o_{\cd\wedge t}) - u^h(t_j, \o_{\cd\wedge t})| &\le& C\sqrt{t_{j}-t_i} \le C\sqrt{t'-t +h}.
\eeaa
These lead to \reff{regt} immediately.
\qed

Combine Lemmas \ref{lem-consistency}, \ref{lem-mon}, and \reff{lem-stability}, it follows from Theorem \ref{thm-PPDE} that
\begin{thm}
\label{thm-conv}
Assume all the conditions in Lemma \ref{lem-stability} hold. Then $u^h$ converges  locally uniformly to the unique viscosity solution $u$ of PPDE \reff{PPDE}. Moreover, 
\bea
\label{reg}
|u(t,\o) - u(t', \o')| \le C\dbf((t,\o), (t',\o')),\q \mbox{for all}~(t,\o), (t',\o')\in\L.
\eea
\end{thm}

\section{The case with classical solution}
\label{sect-rate}
\setcounter{equation}{0}
In this section, we obtain the rate of convergence of our scheme, provided that the PPDE has smooth enough solution. 
Denote
\bea
\label{C12}
C^{2,4}_b:= \Big\{u \in C^{1,2}_b: \pa_t u, \pa_\o u, \pa^2_{\o\o} u \in C^{1,2}_b(\L)\Big\}.
\eea
We shall remark though, as we see in Buckdahn, Ma and Zhang \cite{BMZ},  in general $\pa_t, \pa_{\o_i}, \pa_{\o_j}$ do not commute, and $\pa^2_{\o_i\o_j} u = {1\over 2}[\pa_{\o_i}(\pa_{\o_j} u) + \pa_{\o_j}(\pa_{\o_i}u)]$.

We first have the following general result, in the spirit of Theorem \ref{thm-PPDE}.

\begin{thm}
\label{thm-rate}
Let Assumption \ref{assum-PPDE} hold and  the PPDE \reff{PPDE}  has a classical solution $u\in C^{2,4}_b(\L)$. Assume  a discretization scheme $\dbT^{t,\o}_h$ satisfies: 

(i) For any $(t,\o)\in [0, T)\times \O$ and $\f\in C^{2,4}_b(\L^t)$, 
\bea
\label{consistencyRate}
\Big| {\f(t,\o)-\dbT^{t,\o}_h[\f(t+h,\cd)]\over h}-\cL \f(t,\o)\Big|\leq C h,\q \forall h \in (0, T-t).
\eea

(ii) There exists $L\ge L_0$ such that, for any $(t,\o)\in [0, T)\times \O$  and $\f,\psi\in UC_b(\L^t)$ uniformly Lipschitz continuous,
\bea
\label{monRate}
\Big|\dbT_h^{t,\o}[\f]-\dbT_h^{t,\o}[\psi]\Big|\le (1+Ch)\overline{\cE}_t^{L}[|\f-\psi|].
\eea

\noindent Then we have
\bea
\label{rate}
\max_{0\le i\le n}|u^h(t_i, \o) - u(t_i,\o)| \le Ch.
\eea
\end{thm}
\proof Denote $\e_i := \sup_{\o\in\O} |u^h(t_i,\o) - u(t_i,\o)|$.  Since $\cL u =0$, then \reff{consistencyRate} implies
\beaa
\Big|\dbT^{t_i,\o}_h[u(t_{i+1},\cd)] - u(t_i,\o)\Big|  \le Ch^2.
\eeaa
Now by  \reff{PPDE-uh} and  \reff{monRate} we have
\beaa
|[u^h - u](t_i,\o)| \le \big|\dbT^{t_i,\o}_h[u^h(t_{i+1},\cd)]  - \dbT^{t_i,\o}_h[u(t_{i+1},\cd)] \big| + Ch^{2}\le (1+Ch) \e_{i+1}  + Ch^{2}.
\eeaa
This implies
\bea
\label{ei}
\e_{i} \le (1+Ch) \e_{i+1} + Ch^{2}.
\eea
Since $\e_n = 0$, then by discrete Gronwall inequality we obtain \reff{rate} immediately.
\qed

We next apply the above result to the scheme proposed in Section \ref{sect-scheme}.
\begin{thm}
\label{thm-rate2}
Assume all the conditions in Lemma  \ref{lem-stability} hold true, and the PPDE \reff{PPDE} has a classical solution $u \in C^{2,4}_b(\L)$. Then for the scheme introduced in Section \ref{sect-scheme}, it holds that $|u^h(t_i,\o)-u(t_i, \o)|\le Ch$.
\end{thm}
\proof First, \reff{monRate} follows directly from \reff{Delta2} and Remark \ref{rem-hatP}. By Theorem \ref{thm-rate} it suffices to check \reff{consistencyRate}.  Without loss of generality, we assume $(t,\o) = (0, {\bf 0})$.

 Fix $\f\in C^{2,4}_b(\L)$, and set $\psi := \f(h,\cd)$. Recall the computation in Lemma \ref{lem-consistency} with $(t,\o)=(0,{\bf 0})$ and  $c=0$, we have
\beaa
&& \cD^{(0)} \psi = \f(t+h,{\bf 0}) = \f(0, {\bf 0}) + \int_0^{h} \pa_t \f(s, {\bf 0})ds  \\
&&\q = \f(0, {\bf 0}) + \pa_t\f(0, {\bf 0})  h +  \int_0^{h} \int_0^s \pa_t \pa_t \f(r, {\bf 0})drds =   \f(0, {\bf 0}) + \pa_t\f(0, {\bf 0})  h   + O(h^2);\\
&&{ \cD^{(0)} \psi- \f(0, {\bf 0}) \over h} =    \pa_t\f(0, {\bf 0})     + O(h); \\
&&\cD^{(1)}_i \psi = {1\over \mu_i h} \int_0^{h} \dbE^{\dbP^i} \Big[\big(\pa_t  + \mu_i \pa_{\o^i} \big)\f(s,  B) - \pa_t \f(s, {\bf 0})\Big] ds\\
&&\q= \pa_{\o^i}\f(0,{\bf 0}) + {1\over \mu_i h} \int_0^{h}  \dbE^{\dbP^i} \Big[\big(\pa_t + \mu_i \pa_{\o^i} \big)\f(s,  B)  - \big(\pa_t + \mu_i \pa_{\o^i}\big)\f(0, {\bf 0}) \\
&&\q\qq- [\pa_t \f(s, {\bf 0}) - \pa_t \f(0, {\bf 0})]\Big] ds\\
&&\q= \pa_{\o^i}\f(0,{\bf 0}) + {1\over \mu_i h} \int_0^{h} \int_0^s  \dbE^{\dbP^i} \Big[[\pa_t + \mu_i \pa_{\o^i}]\big(\pa_t  + \mu_i \pa_{\o^i} \big)\f (r,  B) - \pa_t \pa_t \f(r,{\bf 0})\Big] drds\\
&&\q =\pa_{\o^i}\f(0,{\bf 0})  + O(h).
\eeaa
Similarly, we can show that
\beaa
\cD^{(2)}_{i, j} \psi = \pa^2_{\o^i\o^j}\f(0,{\bf 0})+ O(h),\q i,j=1,\cds, d. 
\eeaa
Plug all these into \reff{dbTf} and recall that $G$ is uniformly Lipschitz continuous in $(y,z,\g)$, we obtain \reff{consistencyRate}, and hence prove the theorem.
\qed

\end{document}